\documentclass[a4paper,twoside,11pt]{article}
\usepackage{bbm}
\usepackage{amsfonts}
\usepackage{amssymb}\usepackage{amsmath}
 \newtheorem{rem}{Remark}[section]
 
 \numberwithin{equation}{section}
\begin{document}

\title{Global weak solutions for a periodic two-component $\mu$-Hunter-Saxton system}

\author{Jingjing Liu\footnote{e-mail:
jingjing830306@163.com }\\
Department of Mathematics,
Sun Yat-sen University,\\
510275 Guangzhou, China
\bigskip\\
Zhaoyang Yin\footnote{e-mail: mcsyzy@mail.sysu.edu.cn}
\\ Department of Mathematics, Sun Yat-sen University,\\ 510275 Guangzhou, China}
\date{}
\maketitle

\begin{abstract}
This paper is concerned with global existence of weak solutions
 for a periodic
two-component $\mu$-Hunter-Saxton system. We first derive global
existence for strong solutions to the  system with smooth
approximate initial data.
Then, we show that the limit of approximate solutions is a global weak solution of the two-component $\mu$-Hunter-Saxton system.\\

\noindent 2000 Mathematics Subject Classification: 35G25, 35L05
\smallskip\par
\noindent \textit{Keywords}: A periodic two-component
$\mu$-Hunter-Saxton system,  weak solutions, global existence,
approximate solutions.
\end{abstract}

\section{Introduction}
\newtheorem {definition1}{Definition}[section]
\newtheorem{theorem1}{Theorem}[section]

\par
Recently, a new 2-component system  was introduced by Zuo in
\cite{d-z} as follows:
\begin{equation}
\left\{\begin{array}{ll}
\mu(u)_{t}-u_{txx}=2\mu(u)u_{x}-2u_{x}u_{xx}-uu_{xxx}+\rho\rho_{x}&-\gamma u_{xxx},\\
&t > 0,\,x\in \mathbb{R},\\
\rho_{t}=(\rho u)_x+\gamma \rho_{x}, &t > 0,\,x\in \mathbb{R},\\
u(0,x) = u_{0}(x),& x\in \mathbb{R}, \\
\rho(0,x) = \rho_{0}(x),&x\in \mathbb{R},\\
u(t,x+1)=u(t,x), & t \geq 0, x\in \mathbb{R},\\
\rho(t,x+1)=\rho(t,x), & t \geq 0, x\in \mathbb{R},\\ \end{array}\right. \\
\end{equation}
where $\mu(u)=\int_{\mathbb{S}}udx$ with
$\mathbb{S}=\mathbb{R}/\mathbb{Z}$ and $\gamma \in \mathbb{R}.$ By
integrating both sides of the first equation in the system (1.1)
over the circle $\mathbb{S}=\mathbb{R}/\mathbb{Z}$ and using the
periodicity of $u$, one obtain
$$\mu(u_{t})=\mu(u)_{t}=0.$$
This yields the following periodic 2-component $\mu$-Hunter-Saxton
system:
\begin{equation}
\left\{\begin{array}{ll}
-u_{txx}=2\mu(u)u_{x}-2u_{x}u_{xx}-uu_{xxx}+\rho\rho_{x}&-\gamma u_{xxx},\\
&t > 0,\,x\in \mathbb{R},\\
\rho_{t}=(\rho u)_x+\gamma \rho_{x}, &t > 0,\,x\in \mathbb{R},\\
u(0,x) = u_{0}(x),& x\in \mathbb{R}, \\
\rho(0,x) = \rho_{0}(x),&x\in \mathbb{R},\\
u(t,x+1)=u(t,x), & t \geq 0, x\in \mathbb{R},\\
\rho(t,x+1)=\rho(t,x), & t \geq 0, x\in \mathbb{R},\\ \end{array}\right. \\
\end{equation}
with $\gamma \in \mathbb{R}.$  This system is a 2-component
generalization of the generalized Hunter-Saxton equation obtained in
\cite{k-l-m}. The author \cite{d-z} shows that this system is a
bihamiltonian Euler equation, and also can be viewed as a
bivariational equation.

Obviously, (1.1) is equivalent to (1.2) under the condition
$\mu(u_{t})=\mu(u)_{t}=0.$ In this paper, we will study the system
(1.2) under the assumption $\mu(u_{t})=\mu(u)_{t}=0$.

For $\rho\equiv 0$ and $\gamma=0,$ and replacing $t$ by $-t,$ the
system (1.2) reduces to the generalized Hunter-Saxton equation
(named $\mu$-Hunter-Saxton equation) as follows \cite{k-l-m}:
\begin{equation}
-u_{txx}=-2\mu(u)u_{x}+2u_{x}u_{xx}+uu_{xxx},
\end{equation}
where $\mu(u)=\int_{\mathbb{S}}udx$ denotes the mean of $u.$ The
$\mu$-Hunter-Saxton equation lies mid-way between the periodic
Hunter-Saxton and Camassa-Holm equations with $u=u(t,x)$ being a
time-dependent function on the circle
$\mathbb{S}=\mathbb{R}/\mathbb{Z}$. Recently, the periodic
$\mu$-Hunter-Saxton equation and the periodic
$\mu$-Degasperis-Procesi equation have also been studied in
\cite{escher, fu, l-m-t}. For $\mu(u)=0,$ the equation (1.3) reduces
to the Hunter-Saxton equation \cite{J-R}
\begin{equation}
u_{txx}+2u_{x}u_{xx}+uu_{xxx}=0,
\end{equation}
modeling the propagation of weakly nonlinear orientation waves in a
massive nematic liquid crystal. In the Hunter-Saxton equation
\cite{J-R}, $x$ is the space variable in a reference frame moving
with the linearized wave velocity, $t$ is a slow-time variable and
$u(t,x)$ is a measure of the average orientation of the medium
locally around $x$ at time $t$. More precisely, the orientation of
the molecules is described by the field of unit vectors $(\cos
u(t,x), \sin u(t,x))$ \cite{Y}. The single-component model also
arises in a different physical context as the high-frequency limit
\cite{HHD, J-Y} of the Camassa-Holm equation for shallow water waves
\cite{R-D, J-H} and a re-expression of the geodesic flow on the
diffeomorphism group of the circle \cite{C-K} with a bi-Hamiltonian
structure \cite{F-F} which is completely integrable \cite{C-M}. The
Hunter-Saxton equation also has a bi-Hamiltonian structure
\cite{J-H, P-P} and is completely integrable \cite{B1, J-Y}. The
initial value problem for the Hunter-Saxton equation (1.4) on the
line (nonperiodic case) and on the unit circle
$\mathbb{S}=\mathbb{R}/\mathbb{Z}$ were studied by Hunter and Saxton
in \cite{J-R} using the method of characteristics and by Yin in
\cite{Y} using Kato semigroup method, respectively. Moreover, global
dissipative and conservative weak solutions for the initial boundary
value problem of the Hunter-Saxton equation on the half line were
investigated extensively, c.f. \cite{H-Z1,H-Z2,Z-Z1,Z-Z2,Z-Z3}.

For $\rho\not\equiv0$ , $\gamma=0,$ $\mu(u)=0$ and replacing $t$ by
$-t,$ the system (1.2) reduces to the two-component periodic
Hunter-Saxton system, peakon solutions of the Cauchy problem of it
have been analysed in \cite{C-I}. Moreover, the Cauchy problem and
global weak solutions of two-component periodic Hunter-Saxton system
have been discussed in \cite{l-y} and \cite{g-y} respectively.

Furthermore, the Cauchy problem of (1.2) has been discussed in
\cite{l1-y} recently. The authors established the local
well-posedness, derived precise blow-up scenarios for the system
(1.2) and proved that the system (1.2) has global strong solutions
and also finite time blow-up solutions. However, the existence of
global weak solutions to the system (1.2) has not been studied yet.
The aim of this paper is to present a  global existence result of
weak solutions to the system (1.2).

The main result of this paper is to give the existence of a
global-in-time weak solution $z=\left(
                                                     \begin{array}{c}
                                                       u \\
                                                       \rho \\
                                                     \end{array}
                                                   \right)$ to the problem (1.2) with the initial $z_{0}=\left(
                                                     \begin{array}{c}
                                                       u_{0} \\
                                                       \rho_{0} \\
                                                     \end{array}
                                                   \right)
\in W^{1,\infty}(\mathbb{S})\times L^{\infty}(\mathbb{S}).$  Before
giving the precise statement of the main result, we first introduce
the definition of weak solution to the problem (1.2).

\begin{definition1} $z=\left(
                                                     \begin{array}{c}
                                                       u \\
                                                       \rho \\
                                                     \end{array}
                                                   \right)
\in H^1(\mathbb{S})\times L^{2}(\mathbb{S})$ is said to be an
admissible weak solution to the problem (1.2)  if
$$ z(x,t)\in L_{loc}^{\infty}((0,\infty); H^{1}(\mathbb{S})\times
L^{2}(\mathbb{S}))$$ satisfies the system (1.2) and
$z(t,\cdot)\rightarrow z_0$ as $t\rightarrow 0^+$ in the sense of
distributions on $\mathbb{R}_+\times\mathbb{R}$. Moreover,
$$\|u_x(t,\cdot)\|_{L^2(\mathbb{S})}+\|\rho(t,\cdot)\|_{L^2(\mathbb{S})}\leq
\|u_{0,x}\|_{L^2(\mathbb{S})}+\|\rho_0\|_{L^2(\mathbb{S})}.$$
\end{definition1}
\par
The main result of this paper can be stated as follows:
\begin{theorem1} Let $z_{0}=\left(
                                                     \begin{array}{c}
                                                       u_{0} \\
                                                      \rho_{0} \\
                                                     \end{array}
                                                   \right)
\in W^{1,\infty}(\mathbb{S})\times L^{\infty}(\mathbb{S}).$ If there
exists $\alpha>0$ such that $\rho_0(x)\geq \alpha$ for a.e.
$x\in\mathbb{S},$ then the system (1.2) has an admissible weak
solution
$$ z=\left(
                                                     \begin{array}{c}
                                                       u \\
                                                       \rho \\
                                                     \end{array}
                                                   \right)\in
                                                   C(\mathbb{R}_+;H^1(\mathbb{S})\times L^2(\mathbb{S}))
                                                   \cap L^{\infty}(\mathbb{R}_+;
H^{1}(\mathbb{S})\times L^{2}(\mathbb{S}))$$ in the sense of
Definition 1.1. Furthermore,
$$u\in L^{\infty}_{loc}(\mathbb{R}_+; W^{1,\infty}({\mathbb{S}})) \ \
\ \text{and} \ \ \ \rho\in L^{\infty}_{loc}(\mathbb{R}_+;
L^{\infty}(\mathbb{S})).$$
 \end{theorem1}

\begin{rem}
If there exists $\alpha<0$ such that $\rho_0(x)\leq \alpha$ for a.e.
$x\in\mathbb{S},$ then the conclusions in Theorem 1.1 also hold.
\end{rem}

 The paper is organized as follows. In
 Section 2, we recall some useful lemmas and derive some priori estimates on global strong solutions to (1.2). In
Section 3, we obtain the global existence of approximate solution to
(1.2) with smooth approximate initial data. In Section 4,  acquiring
the precompactness of approximate solutions,
 we prove the
existence of the global weak solution to (1.2).

\section{Preliminaries}
\newtheorem {remark2}{Remark}[section]
\newtheorem{theorem2}{Theorem}[section]
\newtheorem{lemma2}{Lemma}[section]
Since (1.1) is equivalent to (1.2) under the condition
$\mu(u)_{t}=\mu(u_{t})=0,$ to obtain the existence of global weak
solution to (1.2) we study (1.1) henceforth. Moreover, for the sake
of convenience, we let
$$\mu_{0}=\mu(u_{0})=\mu(u)=\int_{\mathbb{S}}u(t,x)dx.$$

We now provide the framework in which we shall reformulate the
system (1.1). We rewrite the system (1.1) as follows:
\begin{equation}
\left\{\begin{array}{ll}
u_{t}-(u+\gamma)u_{x}=\partial_{x}A^{-1}(2\mu_{0}u+\frac{1}{2}u_{x}^{2}+\frac{1}{2}\rho^{2}),
&t > 0,\,x\in \mathbb{R},\\
\rho_{t}-(\rho u)_x=\gamma \rho_{x}, &t > 0,\,x\in \mathbb{R},\\
u(0,x) = u_{0}(x),& x\in \mathbb{R}, \\
\rho(0,x) = \rho_{0}(x),&x\in \mathbb{R},\\
u(t,x+1)=u(t,x), & t \geq 0, x\in \mathbb{R},\\
\rho(t,x+1)=\rho(t,x), & t \geq 0, x\in \mathbb{R},\\ \end{array}\right. \\
\end{equation}
where $A=\mu-\partial_{x}^{2}$ is an isomorphism between
$H^{s}(\mathbb{S})$ and $H^{s-2}(\mathbb{S})$ with the inverse
$v=A^{-1}w$ given explicitly by
\begin{align}
v(x)&=(\frac{x^{2}}{2}-\frac{x}{2}+\frac{13}{12})\mu(w)+(x-\frac{1}{2})\int_{0}^{1}\int_{0}^{y}w(s)dsdy\\
\nonumber&-\int_{0}^{x}\int_{0}^{y}w(s)dsdy+\int_{0}^{1}\int_{0}^{y}\int_{0}^{s}w(r)drdsdy,
\end{align}
which can be found in \cite{escher}. Since $A^{-1}$ and
$\partial_{x}$ commute, the following identities hold
\begin{equation}
A^{-1}\partial_{x}w(x)=(x-\frac{1}{2})\int_{0}^{1}w(x)dx-\int_{0}^{x}w(y)dy+\int_{0}^{1}\int_{0}^{x}w(y)dydx,
\end{equation}
and
\begin{equation}
A^{-1}\partial_{x}^{2}w(x)=-w(x)+\int_{0}^{1}w(x)dx.
\end{equation}
If we rewrite the inverse of the operator $A=\mu-\partial_{x}^{2}$
in terms of a Green's function, we find
$(A^{-1}m)(x)=\int_{0}^{1}g(x-x^{'})m(x^{'})dx^{'}=(g*m)(x).$ So, we
get another equivalent form:
\begin{equation}
\left\{\begin{array}{ll}
u_{t}-(u+\gamma)u_{x}=\partial_{x}g*(2\mu_{0}u+\frac{1}{2}u_{x}^{2}+\frac{1}{2}\rho^{2}),
&t > 0,\,x\in \mathbb{R},\\
\rho_{t}-(\rho u)_x=\gamma \rho_{x}, &t > 0,\,x\in \mathbb{R},\\
u(0,x) = u_{0}(x),& x\in \mathbb{R}, \\
\rho(0,x) = \rho_{0}(x),&x\in \mathbb{R},\\
u(t,x+1)=u(t,x), & t \geq 0, x\in \mathbb{R},\\
\rho(t,x+1)=\rho(t,x), & t \geq 0, x\in \mathbb{R},\\ \end{array}\right. \\
\end{equation}
where the Green's function $g(x)$ is given \cite{l-m-t} by
\begin{equation}
g(x)=\frac{1}{2}x(x-1)+\frac{13}{12} \ \ \ \text{for} \
x\in[0,1)\simeq S^{1},
\end{equation}
and is extended periodically to the real line. In other words,
$$g(x-x^{'})=\frac{(x-x^{'})^{2}}{2}-\frac{|x-x^{'}|}{2}+\frac{13}{12},
\ \ \ x,x^{'}\in[0,1)\simeq S^{1}.$$ In particular, $\mu(g)=1.$
\begin{lemma2}(\cite{constantin, fu}) If $f\in H^{1}(\mathbb{S})$ is such that
$\int_{\mathbb{S}}f(x)dx=0,$ then we have
$$\max\limits_{x\in\mathbb{S}}f^{2}(x)\leq
\frac{1}{12}\int_{\mathbb{S}}f_{x}^{2}(x)dx.$$
\end{lemma2}

Assume that $z=\left(
                                    \begin{array}{c}
                                      u \\
                                      \rho \\
                                    \end{array}
                                  \right)$ is a smooth solution of (2.1). For
convenience, we let
$\mu_{1}=\left(\int_{\mathbb{S}}(u_{0,x}^{2}+\rho_{0}^{2})dx\right)^{\frac{1}{2}}.$
Using the system (1.1), a simple calculation implies
$\frac{d}{dt}\int_{\mathbb{S}}(u_{x}^{2}+\rho^{2})dx=0.$  So we have
\begin{equation}
\mu_{1}=\left(\int_{\mathbb{S}}(u_{0,x}^{2}+\rho_{0}^{2})dx\right)^{\frac{1}{2}}
=\left(\int_{\mathbb{S}}(u_{x}^{2}+\rho^{2})dx\right)^{\frac{1}{2}}.
\end{equation}
Note that $\int_{\mathbb{S}}(u(t,x)-\mu_{0})dx=\mu_{0}-\mu_{0}=0.$
By Lemma 2.1, we find that
$$\max\limits_{x\in\mathbb{S}}[u(t,x)-\mu_{0}]^{2}\leq
 \frac{1}{12}\int_{\mathbb{S}}u_{x}^{2}(t,x)dx\leq
 \frac{1}{12}\mu_{1}^{2}.$$ So we have
 \begin{equation}
 \|u(t,\cdot)\|_{L^{\infty}(\mathbb{S})}\leq
 |\mu_{0}|+\frac{\sqrt{3}}{6}\mu_{1}.
 \end{equation}
\begin{lemma2} (\cite{l1-y}) Let $z_{0}=\left(
                                                     \begin{array}{c}
                                                       u_{0} \\
                                                       \rho_{0} \\
                                                     \end{array}
                                                   \right)
\in H^2(\mathbb{S})\times H^{1}(\mathbb{S}).$ If $\rho_0(x)\neq0 $
 for all $x\in\mathbb{S}$, then the corresponding strong solution $z=\left(
                                    \begin{array}{c}
                                      u \\
                                      \rho \\
                                    \end{array}
                                  \right)$
 to (1.1) exists globally in time, i.e. $z\in
 C(\mathbb{R}_+;H^2(\mathbb{S})\times H^{1}(\mathbb{S}))\cap C^1(\mathbb{R}_+;H^{1}(\mathbb{S})\times L^{2}(\mathbb{S}))$.
 Moreover,
 there exists $\beta>0$ such that for all $t\in\mathbb{R}_{+},$
\begin{align*}
\|u_x(t,\cdot)\|_{L^{\infty}(\mathbb{S})}\leq\frac{1}{2\beta}(1+\|
\rho_0\|^2_{L^{\infty}(\mathbb{S})}&+\|
u_{0,x}\|^2_{L^{\infty}(\mathbb{S})})\\
&\cdot e^{(4\mu_{0}^{2}+\frac{1}{2}\mu_{1}^{2}+\frac{\sqrt{3}}{3}
|\mu_{0}|\mu_{1}+\frac{1}{2})t}:= C_1(t)
\end{align*}
 and
$$\|\rho(t,\cdot)\|_{L^{\infty}(\mathbb{S})}\leq e^{C_1(t)t}\|\rho_0\|_{L^{\infty}(\mathbb{S})}:= C_2(t),$$
where $\beta= \inf_{x\in\mathbb{S}}|\rho_0(x)|.$
\end{lemma2}

\begin{lemma2}(\cite{s}) Assume $X\subset B \subset Y$ with compact
imbedding $X\rightarrow B$ ($X,B$ and $Y$ are Banach spaces), $1\leq
p\leq \infty$ and (1) $F$ is bounded in $L^{p}(0,T; X)$, (2)
$\|\tau_{h}f-f\|_{L^{p}(0,T-h; Y)}\rightarrow 0$ as $h\rightarrow 0$
uniformly for $f\in F.$ Then $F$ is relatively compact in
$L^{p}(0,T; B)$ (and in $C(0,T; B)$ if $p=\infty$), where
$(\tau_{h}f)(t)=f(t+h)$ for $h>0,$ if $f$ is defined on $[0,T],$
then the translated function $\tau_{h}f$ is defined on $[-h, T-h].$
\end{lemma2}

\begin{lemma2}(Appendix C of \cite{Li}) Let $X$ be a separable
reflexive Banach space and let $f^n$ be bounded in
$L^{\infty}(0,T;X)$ for some $T\in(0,\infty).$ We assume that
$f^n\in C([0,T];Y)$ where $Y$ is a Banach space such that
$X\hookrightarrow Y$, $Y'$ is separable and dense in $X'$.
Furthermore, $(\phi,f^n(t))_{Y'\times Y}$ is uniformly continuous in
$t\in [0,T]$ and uniformly in $n\geq 1.$ Then $f^n$ is relatively
compact in $C^{w}([0,T];X)$, the space of continuous functions from
$[0,T]$ with values in $X$ when the latter space is equipped with
its weak topology.
\end{lemma2}

\begin{remark2} If the conditions which $f^n$ satisfies in Lemma 2.4
are replaced by the following conditions:
$$f^n\in L^{\infty}(0,T;X),\ \partial_t f^n\in L^p(0,T;Y) \text{  for
some}\,\, p\in(1,\infty],$$ and $$ \| f^n\|_{L^{\infty}(0,T;X)}, \
\|
\partial_tf^n\|_{L^p(0,T;Y)}\leq C, \ \ \forall n\geq1,$$
then the conclusion of Lemma 2.4 holds true.
\end{remark2}

\section{Global approximate solutions}
\newtheorem {remark3}{Remark}[section]
\newtheorem{theorem3}{Theorem}[section]
\newtheorem{lemma3}{Lemma}[section]
\newtheorem{claim3}{Claim}[section]

In the section, we first prove the existence of approximate
solutions. Then, with the basic estimates given in Section 2, we
will give some useful estimates to the approximate solutions.

\par
Let $z_{0}=\left(
                                                     \begin{array}{c}
                                                       u_{0} \\
                                                      \rho_{0} \\
                                                     \end{array}
                                                   \right)
\in W^{1,\infty}(\mathbb{S})\times L^{\infty}(\mathbb{S})$ and there
exists $\alpha>0$ such that $\rho_0(x)\geq \alpha$ for a.e.
$x\in\mathbb{S}.$ Define $z_0^n:=\left(
                                                     \begin{array}{c}
                                                       \phi_n\ast u_{0} \\
                                                      \phi_n\ast \rho_{0} \\
                                                     \end{array}
                                                   \right)= \left(
                                                     \begin{array}{c}
                                                        u_{0}^n \\
                                                      \rho_{0}^n \\
                                                     \end{array}
                                                   \right)\in H^2(\mathbb{S})\times H^{1}(\mathbb{S})$, for $n\geq
                                                   1$,
here $\{\phi_n\}_{n\geq 1}$ are the mollifiers
$$ \phi_n(x):=\left(\int_{\mathbb{R}}\phi(\xi)d\xi\right)^{-1}n\phi(nx),\ \ \ \ \ x\in\mathbb{R}, \ \ n\geq 1,$$
where $\phi\in C_c^{\infty}(\mathbb{R})$ is defined by
$$\phi(x)=\left \{\begin {array}{ll}e^{1/(x^2-1)},  \ \ \ \ \ \ \ \ |x|<1,
 \\ 0, \ \ \ \ \ \ \ \ \ \ \ \ \ \ \ \ \ \ |x|\geq 1.\end {array}\right.$$

In view of $\rho_{0}(x)\geq\alpha>0,$ for a.e. $x\in\mathbb{S}$ and
$\phi_n(x)\geq0,$ we have  $$ \rho_{0}^n(x)=\phi_n\ast
\rho_{0}(x)\geq \alpha \int_{\mathbb{R}}\phi_n(y)dy=\alpha>0, \ \
\forall x\in \mathbb{S}.$$ Clearly, we also have
\begin{align}
 u_0^n\rightarrow u_0 \ \ \text{in}  \ \
H^1(\mathbb{S}), \ \ \ \ \rho_0^n\rightarrow \rho_0 \ \ \text{in}\ \
L^2(\mathbb{S}),\ \ \ \text{as}\ \ n\rightarrow\infty
\end{align} and
\begin{align}
\|u_0^n\|_{L^2(\mathbb{S})}\leq \|u_0\|_{ L^2(\mathbb{S})},
\|u_{0,x}^n\|_{L^2(\mathbb{S})}\leq \|u_{0,x}\|_{L^2(\mathbb{S})},
\|\rho_0^n\|_{L^2(\mathbb{S})}\leq \|\rho_0\|_{L^2(\mathbb{S})}
.\end{align} Thus, we obtain the corresponding solution $z^n\in
 C(\mathbb{R}_+;H^2(\mathbb{S})\times H^{1}(\mathbb{S}))\cap C^1(\mathbb{R}_+;H^{1}(\mathbb{S})\times L^{2}(\mathbb{S}))$ to the system (1.1)
 with initial data $z^n_0$ by Lemma 2.2 under the condition
 $\mu(u^{n})_{t}=\mu(u^{n}_{t})=0.$

 For given $z_0=\left(
                                                     \begin{array}{c}
                                                       u_{0} \\
                                                       \rho_{0} \\
                                                     \end{array}
                                                   \right)
\in  W^{1,\infty}(\mathbb{S})\times L^{\infty}(\mathbb{S})$, we set
$$\mu^{n}_{0}:=\mu(u^{n}_{0})=\mu(u^{n})=\int_{\mathbb{S}}u^{n}(t,x)dx,$$
$$\mu^{n}_{1}:=\left(\int_{\mathbb{S}}((u^{n}_{0,x})^{2}+(\rho^{n}_{0})^{2})dx\right)^{\frac{1}{2}}
=\left(\int_{\mathbb{S}}((u^{n}_{x})^{2}+(\rho^{n})^{2})dx\right)^{\frac{1}{2}}$$
 Then we have the following remark.

\begin{remark3}
By (2.7)-(2.8) and (3.2), we have
\begin{equation}
\mu^{n}_{0}\rightarrow \mu_{0} \ \text{as}\
 n\rightarrow\infty, \ \ \ \ \ \ \ \ \ (\mu^{n}_{1})^{2}\rightarrow \mu_{1}^{2} \ \text{as}\
 n\rightarrow\infty,
\end{equation}
and
\begin{align}\| u_x^n(t,\cdot)\|_{L^2(\mathbb{S})}^2+\|
\rho^n(t,\cdot)\|^2_{L^2(\mathbb{S})}=(\mu^{n}_{1})^{2}\leq\mu_{1}^{2},
\
 \ \ \forall t\in\mathbb{R}_+.
\end{align}
Moreover, we get
\begin{equation}
|\mu_{0}^{n}|\leq\int_{\mathbb{S}}|u_{0}^{n}|dx\leq\|u_{0}^{n}\|_{L^2(\mathbb{S})}\leq\|u_{0}\|_{L^2(\mathbb{S})},
\end{equation}
\begin{align}
 \|u^n(t,\cdot)\|_{L^\infty(\mathbb{S})}\leq
 |\mu_{0}^{n}|+\frac{\sqrt{3}}{6}\mu_{1}^{n}\leq
 \|u_{0}\|_{L^{2}(\mathbb{S})}+\frac{\sqrt{3}}{6}\mu_{1}.
 \end{align}
Furthermore, by Lemma 2.2, we obtain
\begin{align}
\|u^{n}_x(t,\cdot)\|_{L^{\infty}(\mathbb{S})}\leq&\frac{1}{2\beta}(1+\|
\rho_0\|^2_{L^{\infty}(\mathbb{S})}+\|
u_{0,x}\|^2_{L^{\infty}(\mathbb{S})})\\
\nonumber&\cdot
e^{(4\|u_{0}\|_{L^{2}(\mathbb{S})}^{2}+\frac{1}{2}\mu_{1}^{2}+\frac{\sqrt{3}}{3}
\|u_{0}\|_{L^{2}(\mathbb{S})}\mu_{1}+\frac{1}{2})t}:=
\widetilde{C_1}(t)
\end{align}
 and
\begin{equation}
\|\rho^{n}(t,\cdot)\|_{L^{\infty}(\mathbb{S})}\leq
e^{\widetilde{C_1}(t)t}\|\rho_0\|_{L^{\infty}(\mathbb{S})}:=
\widetilde{C_2}(t),
\end{equation}
where $\beta= \inf_{x\in\mathbb{S}}|\rho_0(x)|$.
\end{remark3}

\section{The existence of global weak solution}
\newtheorem {remark4}{Remark}[section]
\newtheorem{theorem4}{Theorem}[section]
\newtheorem{lemma4}{Lemma}[section]
\newtheorem{claim4}{Claim}[section]
In this section, with the basic energy estimate in Section 3, we are
ready to obtain the necessary compactness of the approximate
solutions $z^n(t,x).$ Acquiring the precompactness of approximate
solutions, we prove the existence of the global weak solution to the
system (1.1).
\begin{lemma4}
For any fixed $T>0$, there exist a subsequence $\{z^{n_k}(t,x)\}$ of
the sequence $\{z^{n}(t,x)\}$ and some function  $z(t,x)\in
L^\infty((0,\infty); H^1(\mathbb{S})\times
L^2(\mathbb{S}))\cap(H^{1}((0,T)\times \mathbb{S})\times
L^{2}((0,T)\times \mathbb{S})$ such that
\begin{equation} z^{n_k}\rightharpoonup z \ \ \ \ \text{in}\
H^1((0,T)\times\mathbb{S})\times L^2((0,T)\times\mathbb{S}) \
\text{as}\ n_k\rightarrow\infty, \forall T>0,
\end{equation}
and
\begin{equation} u^{n_k}\rightarrow u \ \ \text{in} \ \
L^\infty((0,T)\times\mathbb{S}) \ \
 \text{as} \ \ n_k\rightarrow\infty.
\end{equation}
Moreover, $u(t,x)\in C((0,T);L^{\infty}(\mathbb{S})).$
\end{lemma4}
\textbf{Proof}
It follows from Remark 3.1 that $\{z^{n}(t,x)\}$ is uniformly
bounded in $L^\infty((0,\infty); H^1(\mathbb{S})\times
L^2(\mathbb{S}))$.

We will prove that the sequence $\{z^{n}(t,x)\}$ is uniformly
bounded in the space $H^{1}((0,T)\times \mathbb{S})\times
L^{2}((0,T)\times \mathbb{S}).$ Firstly, we claim that
$\{u^{n}(t,x)\}$ is uniformly bounded in the space
$H^{1}((0,T)\times \mathbb{S}).$ Using Remark 3.1, we have
\begin{align*}
\|(u^{n}+\gamma)u^{n}_{x}\|_{L^{2}((0,T)\times \mathbb{S})}&\leq
 \|u^{n}u^{n}_{x}\|_{L^{2}((0,T)\times \mathbb{S})}
 +|\gamma|\|u^{n}_{x}\|_{L^{2}((0,T)\times \mathbb{S})}\\
 &\leq(\|u_{0}\|_{L^{2}(\mathbb{S})}+\frac{\sqrt{3}}{6}\mu_{1}+|\gamma|)\mu_{1}\sqrt{T},
\end{align*}
and
\begin{align*}
&\|g_{x}*(2\mu_{0}^{n}u^{n}+\frac{1}{2}(u_{x}^{n})^{2}+\frac{1}{2}(\rho^{n})^{2}\|_{L^{2}((0,T)\times
\mathbb{S})}\\
\leq \ & \|g_{x}\|_{L^{1}((0,T)\times
\mathbb{S})}\|2\mu_{0}^{n}u^{n}+\frac{1}{2}(u_{x}^{n})^{2}
+\frac{1}{2}(\rho^{n})^{2}\|_{L^{2}((0,T)\times
\mathbb{S})}\\
\leq \ &
\frac{T}{2}\left(\int_{0}^{T}\int_{\mathbb{S}}(\mu_{0}^{n})^{2}dxdt
+\int_{0}^{T}\int_{\mathbb{S}}(u^{n})^{2}dxdt+\frac{1}{2}\int_{0}^{T}\int_{\mathbb{S}}
\left((u_{x}^{n})^{2}+(\rho^{n})^{2}\right)dxdt\right)\\
\leq \ &
\frac{T}{2}\left[\|u_{0}\|_{L^{2}(\mathbb{S})}^{2}T+\left(\|u_{0}\|_{L^{2}(\mathbb{S})}
+\frac{\sqrt{3}}{6}\mu_{1}\right)^{2}+\frac{1}{2}\mu_{1}^{2}T\right].
\end{align*}
Then, by the first equation in (2.2), we know that
$\{u^{n}_{t}(t,x)\}$ is uniformly bounded in $L^{2}((0,T)\times
\mathbb{S}).$ Combing this conclusion with Remark 3.1, we obtain
$\{u^{n}(t,x)\}$ is uniformly bounded in the space
$H^{1}((0,T)\times \mathbb{S}).$ Furthermore, from Remark 3.1, one
can easily obtain that $\rho^n$ is uniformly bounded in the space
$L^2((0,T)\times\mathbb{S})$, and thus (4.1) follows.

Observe that, for each $0\leq s, t\leq T,$
$$\|u^{n}(t,\cdot)-u^{n}(s,\cdot)\|_{L^{2}(\mathbb{S})}^{2}
=\int_{\mathbb{S}}(\int_{s}^{t}\frac{\partial u^{n}}{\partial
\tau}(\tau, x)d\tau)^{2}dx\leq
|t-s|\int_{0}^{T}\int_{\mathbb{S}}(u_{t}^{n})^{2}dxdt.$$ Moreover,
$\{u^{n}(t,x)\}$ is uniformly bounded in $L^{\infty}(0,T;
H^{1}(\mathbb{S}))$ and $H^{1}(\mathbb{S})\subset\subset
L^{\infty}(\mathbb{S})\subset L^{2}(\mathbb{S}),$ then (4.2) and
$u(t,x)\in C((0,T);L^{\infty}(\mathbb{S}))$ is consequence of Lemma
2.3.

\begin{remark4} By Remark 3.1 and the above argument,  there exists a subsequences of
$\{(u_x^n)^2,(\rho^n)^2\}$, denoted again by
$\{(u_x^{n_k})^2,(\rho^{n_k})^2\}$, converging weakly in
$L_{loc}^p(\mathbb{R}_+\times\mathbb{R})$ ,  where $1<p<\infty$,
i.e. there exist $\overline{u_x^2}\in
L_{loc}^p(\mathbb{R}_+\times\mathbb{R})$ and $\overline{\rho^2}\in
L_{loc}^p(\mathbb{R}_+\times\mathbb{R})$ such that
$$ (u_x^{n_k})^2\rightharpoonup \overline{u_x^2} \ \ \ \text{and} \ \
(\rho^{n_k})^2\rightharpoonup\overline{\rho^2} \ \ \text{in}\ \
L_{loc}^p(\mathbb{R}_+\times\mathbb{R}).$$ Moreover,
 we have
 \begin{equation*}
u_x^{n_k}\rightharpoonup u_x \ \ \text{in} \
L_{loc}^p(\mathbb{R}_+\times\mathbb{R})\ \text{and}\ \ \
u_x^{n_k}\rightharpoonup^{\star} u_x \ \ \text{in}\ \
L_{loc}^{\infty}(\mathbb{R}_+;L^2( \mathbb{S})),
\end{equation*}
\begin{equation*}
\rho^{n_k}\rightharpoonup \rho \ \ \text{in} \
L_{loc}^p(\mathbb{R}_+\times\mathbb{R})\ \ \text{and}\ \ \
\rho^{n_k}\rightharpoonup^{\star} \rho \ \ \text{in}\ \
L_{loc}^{\infty}(\mathbb{R}_+;L^2( \mathbb{S})),
\end{equation*}
 Furthermore, we have
\begin{equation}
u_x^2(t,x)\leq \overline{u_x^2}(t,x), \ \ \ \
\rho^2(t,x)\leq\overline{\rho^2}(t,x) \ \ a.e. \ \text{on}\
\mathbb{R}_+\times\mathbb{R}.
\end{equation}
In view of (3.7) and (3.8), we have
\begin{align}\|u_x(t,\cdot)\|_{L^\infty(\mathbb{S})}\leq
\widetilde{C_1}(t),\ \   \ \
\|\rho(t,\cdot)\|_{L^\infty(\mathbb{S})}\leq \widetilde{C_2}(t),\ \
\forall t\in\mathbb{R}_+,
\end{align}
where $\widetilde{C_1}(t)$ and $\widetilde{C_2}(t)$ are given in
Remark 3.1.
\end{remark4}
\par

\begin{lemma4} There hold
\begin{align}
\lim_{t\rightarrow0^+}\int_{\mathbb{S}}u_x^2(t,x)dx=
\lim_{t\rightarrow0^+}\int_{\mathbb{S}}\overline{u_x^2}(t,x)=\int_{\mathbb{S}}u_{0,x}^2(x)dx
\end{align}
and
\begin{eqnarray}
\lim_{t\rightarrow0^+}\int_{\mathbb{S}}\rho^2(t,x)dx=
\lim_{t\rightarrow0^+}\int_{\mathbb{S}}\overline{\rho^2}(t,x)dx=\int_{\mathbb{S}}\rho_0^2(x)dx.
\end{eqnarray}
\end{lemma4}
\textbf{Proof}
By Lemma 4.1, for any $T>0$, we have $u^n\in
L^{\infty}((0,T);H^1(\mathbb{S}))$, $u_t^n$ are uniformly bounded in
$L^{\infty}((0,T); L^2(\mathbb{S}))$. By Lemma 2.2, we have $u^n\in
C([0,T];H^1(\mathbb{S}))$. Then in view of Lemma 2.4, Remark 2.1 and
the proof of Lemma 4.1, we have $\{u^{n}\}$ contains a subsequence,
denoted again by $\{u^{n_k}\}$, converging weakly in
$H^1(\mathbb{S})$ uniformly in $t\in(0,T)$. The limit function is
$u$. This implies that $u$ is weakly continuous from $(0,T)$ into
$H^1(\mathbb{S})$, i.e.,
\begin{equation}
u\in C^w([0,T];H^1(\mathbb{S})).
\end{equation}
 Similarly, as
$\rho^{n}\in L^{\infty}((0,T);L^2(\mathbb{S}))$, in view of (3.4)
and (3.6), we get that for all $t\in(0,T)$,
\begin{eqnarray*}
\|\rho^n_t(t,\cdot)\|_{H^{-1}(\mathbb{S})}&=&\sup_{\|
\phi\|_{H^1(\mathbb{S})}=1}\int_{\mathbb{S}}((u^n\rho^n)_{x}\phi+\gamma
\rho^{n}_{x}\phi )dx\leq\|
u^n\rho^n\|_{L^2(\mathbb{S})}+|\gamma|\|\rho^{n}\|_{L^2(\mathbb{S})}\\
&\leq&(\|u_{0}\|_{L^2(\mathbb{S})}+\frac{\sqrt{3}}{6}\mu_{1}+|\gamma|)\mu_{1}.
\end{eqnarray*}
This shows that $\rho_t^n$ is uniformly bounded in $
L^{\infty}((0,T); H^{-1}(\mathbb{S}))$. Then by Lemma 2.4 and Remark
2.1, we have $\{\rho^{n}\}$ contains a subsequence, denoted again by
$\{\rho^{n_k}\},$ converging weakly in $L^2(\mathbb{S})$ uniformly
in $t$. The limit function is $\rho$. This implies that $\rho$ is
weakly continuous from $(0,T)$ into $L^2(\mathbb{S})$, i.e.,
\begin{equation}
\rho\in C^w([0,T];L^2(\mathbb{S})).
\end{equation}
Then by (4.7)-(4.8),  we get
$$ \rho(t,\cdot)\rightharpoonup \rho_0 \ \ \text{and} \
\ u_x(t,\cdot)\rightharpoonup u_{0,x} \ \ \ \text{in}\ \
L^2(\mathbb{S})\ \ \ \text{as}\ \ t\rightarrow0^+.$$ Thus, we have
\begin{eqnarray*}
\liminf_{t\rightarrow0^+}\int_{\mathbb{S}}\rho^2(t,x)dx\geq\int_{\mathbb{S}}\rho_0^2(x)dx
\end{eqnarray*}
and
\begin{eqnarray*}
\liminf_{t\rightarrow0^+}\int_{\mathbb{S}}u_x^2(t,x)dx\geq\int_{\mathbb{S}}u_{0,x}^2(x)dx.
\end{eqnarray*}
Therefore, we deduce
\begin{eqnarray}
&&\liminf_{t\rightarrow0^+}\int_{\mathbb{S}}(u_x^2(t,x)+\rho^2(t,x))dx\\&\geq&
\liminf_{t\rightarrow0^+}\int_{\mathbb{S}}u_x^2(t,x)dx+\liminf_{t\rightarrow0^+}\int_{\mathbb{S}}\rho^2(t,x)dx
\nonumber\\&\geq& \nonumber
\int_{\mathbb{S}}(u_{0,x}^2(x)+\rho_0^2(x))dx.
\end{eqnarray}
Moreover, by Remark 4.1 and (3.2), we infer
\begin{eqnarray*}
&&\int_{\mathbb{S}}\left(\overline{u_x^2}(t,x)+\overline{\rho^2}(t,x)\right)dx\\&\leq&
\liminf_{n_k\rightarrow\infty}\int_{\mathbb{S}}\left((u^{n_k}_x)^2(t,x)+(\rho^{n_k})^2(t,x)\right)dx\\&=&
\liminf_{n_k\rightarrow\infty}\int_{\mathbb{S}}\left((u^{n_k}_{0,x})^2(x)+(\rho^{n_k}_0)^2(x)\right)dx\\&\leq&
\int_{\mathbb{S}}(u_{0,x}^2(x)+\rho_0^2(x))dx.
\end{eqnarray*}
 Thus, we obtain
\begin{align}
\limsup_{t\rightarrow0^+}\int_{\mathbb{S}}\left(\overline{u_x^2}(t,x)+
\overline{\rho^2}(t,x)\right)dx\leq\int_{\mathbb{S}}(u_{0,x}^2(x)+\rho_0^2(x))dx.
\end{align}
In view of (4.3), (4.9)-(4.10), we get
\begin{eqnarray*}
\liminf_{t\rightarrow0^+}\int_{\mathbb{S}}(u_x^2(t,x)+\rho^2(t,x))dx
&\geq&
\int_{\mathbb{S}}(u_{0,x}^2(x)+\rho_0^2(x))dx\\
&\geq&\limsup_{t\rightarrow0^+}\int_{\mathbb{S}}\left(\overline{u_x^2}(t,x)+
\overline{\rho^2}(t,x)\right)dx\\
&\geq&\limsup_{t\rightarrow0^+}\int_{\mathbb{S}}(u_x^2(t,x)+\rho^2(t,x))dx.
\end{eqnarray*}
This completes the proof of the lemma.
\par

\begin{lemma4}There hold
\begin{equation}
\frac{\partial}{\partial
t}(\overline{u_{x}^{2}}+\overline{\rho^{2}})-\frac{\partial}{\partial
x}[(u+\gamma)(\overline{u_{x}^{2}}+\overline{\rho^{2}})]=-4\mu_{0}uu_{x}+4\mu_{0}^{2}u_{x}+u_{x}\mu_{1}^{2}
\end{equation}
in the sense of distributions on $((0,T)\times\mathbb{R}).$
\end{lemma4}
\textbf{Proof} Note that $z^{n_{k}}$ is the solution of the system
(2.1) with the initial $z_{0}^{n_{k}}.$ Differentiating the first
equation in (2.1) with respect $x$ and using
$\partial_{x}^{2}A^{-1}w=-w+\mu(w)$, we have
$$u^{n}_{tx}-(u^{n}+\gamma)u^{n}_{xx}=-2\mu^{n}_{0}u^{n}+\frac{1}{2}(u^{n}_{x})^{2}-\frac{1}{2}(\rho^{n})^{2}
+2(\mu_{0}^{n})^{2}+\frac{1}{2}(\mu_{1}^{n})^{2}.$$ Then
\begin{align*}
\frac{\partial}{\partial t} (u_x^{n_k})^2-\frac{\partial}{\partial
x}((u^{n_k}+\gamma)(u_x^{n_k})^2)=-4\mu_{0}^{n_k}u^{n_k}u_{x}^{n_k}&-u_{x}^{n_k}(\rho^{n_k})^{2}\\
&+4(\mu_{0}^{n_k})^{2}u_{x}^{n_k}+u_{x}^{n_k}(\mu_{1}^{n_k})^{2}
\end{align*}
and
$$\frac{\partial}{\partial t}(\rho^{n_k})^2-
\frac{\partial}{\partial
x}((u^{n_k}+\gamma)(\rho^{n_k})^2)=u_x^{n_k}(\rho^{n_k})^2.$$ Adding
the above two equalities and  letting $n_k\rightarrow\infty$, in
view of Lemma 4.1 and Remark 4.1, we get (4.11).

\begin{lemma4} There hold
\begin{align}
&\frac{\partial}{\partial
t}(u_x^2+\rho^{2})-\frac{\partial}{\partial
x}((u+\gamma)(u_x^2+\rho^{2}))\\
\nonumber= \ &u_{x}^{3}-4\mu_{0}uu_{x}-u_{x}\overline{u_{x}^{2}}
-u_{x}\overline{\rho^{2}}+4\mu_{0}^{2}u_{x}+\mu_{1}^{2}u_{x}+u_{x}\rho^{2}
\end{align}
in the sense of distributions on $((0,T)\times\mathbb{R}).$
\end{lemma4}
\textbf{Proof} Since $z^{n_{k}}$ satisfy
$$u^{n}_{tx}-((u^{n}+\gamma)u^{n}_{x})_{x}=-2\mu^{n}_{0}u^{n}-\frac{1}{2}(u^{n}_{x})^{2}
-\frac{1}{2}(\rho^{n})^{2}+2(\mu^{n}_{0})^{2}+\frac{1}{2}(\mu^{n}_{1})^{2}$$
and
$$\rho^{n}_{t}-(\rho^{n} u^{n})_{x}=\gamma \rho^{n}_{x}.$$ In view of Lemma 4.1 and
Remark 4.1, letting $n_k\rightarrow\infty$, we obtain
\begin{equation}
u_{tx}-((u+\gamma)u_{x})_{x}=-2\mu_{0}u-\frac{1}{2}\overline{u_{x}^{2}}-\frac{1}{2}\overline{\rho^{2}}
+2\mu_{0}^{2}+\frac{1}{2}\mu_{1}^{2}
\end{equation}
 and
\begin{equation}
\rho_{t}-(\rho u)_{x}=\gamma \rho_{x}
\end{equation}
in the sense of distributions on $((0,T)\times\mathbb{R}).$

Denote $u_{n,x}(t,x):=\left((u_x(t,\cdot)\ast \phi_{n}\right)(x)$
and $\rho_{n}(t,x):=(\rho(t,\cdot)\ast \phi_n)(x).$ According to
Lemma II.1 of \cite{D-L}, it follows from (4.13)-(4.14) that
$u_{n,x}$ and $\rho_{n}$ solve
\begin{equation}
\frac{\partial u_{n,x}}{\partial t}-(u+\gamma)\frac{\partial
u_{n,x}}{\partial
x}=(u^2_x-2\mu_{0}u-\frac{1}{2}\overline{u_{x}^{2}}-\frac{1}{2}\overline{\rho^{2}}
+2\mu_{0}^{2}+\frac{1}{2}\mu_{1}^{2})\ast\phi_{n}+\tau_{n}
\end{equation} and
\begin{equation}
\frac{\partial\rho_{n}}{\partial
t}-(u+\gamma)\frac{\partial\rho_{n}}{\partial
x}=(u_x\rho)\ast\phi_{n}+\sigma_{n},
\end{equation}
where the errors $\tau_n=\left(u\frac{\partial u_{x}}{\partial
x}\right)\ast\phi_{n}-u\cdot\frac{\partial u_{n,x}}{\partial x}$ and
$\sigma_n=\left(u\frac{\partial \rho}{\partial
x}\right)\ast\phi_{n}-u\cdot\frac{\partial \rho_{n,x}}{\partial x}$
tend to zero in $L^1_{loc}(\mathbb{R}_+\times\mathbb{R}).$ Using
(4.15) and (4.16), we get
\begin{align}
&\frac{\partial(u_{n,x})^2}{\partial t}-\frac{\partial}{\partial
x}((u+\gamma)u_{n,x}^2)\\
\nonumber= \
&2u_{n,x}\left((u^2_x-2\mu_{0}u-\frac{1}{2}\overline{u_{x}^{2}}-\frac{1}{2}\overline{\rho^{2}}
+2\mu_{0}^{2}+\frac{1}{2}\mu_{1}^{2})\ast\phi_{n}+\tau_{n}\right)-u_{x}u_{n,x}^{2}
\end{align}
and
\begin{align}
\frac{\partial(\rho_{n})^2}{\partial t}-\frac{\partial}{\partial
x}((u+\gamma)\rho_{n}^2)=2\rho_{n}((u_x\rho)\ast\phi_{n}+\sigma_{n})-u_{x}\rho_{n}^{2}.
\end{align}
Sending $n\rightarrow\infty$ in (4.17) and (4.18) and adding the
results yield (4.12).

Next, we give the main lemma of this section.
\begin{lemma4}
There hold
\begin{equation}
\overline{u_x^2}(t,x)=u_x^2(t,x) \ \ \ \text{and} \ \ \
\overline{\rho^2}(t,x)=\rho^2(t,x) \ \ a.e. \ \text{on}\ \
(0,T)\times\mathbb{R}.
\end{equation}
\end{lemma4}
\textbf{Proof}
Subtracting (4.12) from (4.11) and  integrating the obtained
equality over $(\varepsilon,t)\times\mathbb{S})$ give
\begin{align*}
&\int_{\mathbb{S}}\left(\overline{u_x^2}-u_x^2+\overline{\rho^2}-\rho^2\right)(t,x)dx\\
=&
\int_{\varepsilon}^t\int_{\mathbb{S}}\left(\overline{u_x^2}-u_x^2+\overline{\rho^2}
-\rho^2\right)u_x(t,x)dxds\\
&+\int_{\mathbb{S}}\left(\overline{u_x^2}(\varepsilon,x)-u_x^2(\varepsilon,x)+\overline{\rho^2}(\varepsilon,x)-\rho^2(\varepsilon,x)\right)dx,
\end{align*}
for almost all $t\in(0,T)$. Letting $\varepsilon\rightarrow0$ and
using Lemma 4.2 and (4.4), we get
\begin{align*}
\int_{\mathbb{S}}\left(\overline{u_x^2}-u_x^2+\overline{\rho^2}-\rho^2\right)(t,x)dx
\leq
\widetilde{C_1}(T)\int_{0}^t\int_{\mathbb{S}}\left(\overline{u_x^2}-u_x^2+\overline{\rho^2}-\rho^2\right)dxds.
\end{align*}
Using Gronwall's inequality,  we obtain
\begin{align*}\int_{\mathbb{S}}\left(\overline{u_x^2}-u_x^2+\overline{\rho^2}-\bar{\rho}^2\right)(t,x)dx
\leq0.
\end{align*}
By (4.3), we deduce
\begin{align*}
0\leq\int_{\mathbb{S}}\left(\overline{u_x^2}-u_x^2+\overline{\rho^2}-\rho^2\right)(t,x)dx\leq0,
\end{align*}
which yields
$$\int_{\mathbb{S}}(\overline{u_x^2}-u_x^2)(t,x)dx=\int_{\mathbb{S}}(\overline{\rho^2}-\rho^2)(t,x)dx=0.$$
In view of $u^n_x, u_x$ and $\rho^n,\rho$ being periodic with
respect to $x$, we deduce that (4.19) holds.

To prove our main theorem, we need the following lemma.

\begin{lemma4} (\cite{E-L}) Let U be the bounded domain in $\mathbb{R}^n.$
Assume that the sequence $\{f_k\}_{k=1}^{\infty}$ is bounded in
$L^\infty(U;\mathbb{R}^m)$. Then there exist a subsequence
$\{f_{k_j}\}_{j=1}^{\infty}\subset\{f_k\}_{k=1}^{\infty}$ and for
a.e. $x\in U,$ a Borel probability measure $\nu_x$ on $\mathbb{R}^m$
such that for $F\in C(\mathbb{R}^m)$ we have
$$F(f_{k_j})\rightharpoonup^*\overline{F} \ \ \ \ \text{in}\ \
L^\infty(U),$$ where $\overline{F}(x)=\int_{\mathbb{R}^m}F(y)d\nu_x$
for a.e. $x\in U$. Moreover, if $\nu_x$ is a unit point mass for
a.e. $x\in U$, then $f_{k_j}\rightarrow f$ in $L^2(U;\mathbb{R}^m)$,
here $f(x)=\int_{\mathbb{R}^m}yd\nu_x$.
\end{lemma4}

By (3.7)-(3.8), for any $T>0$ and $a<b$, $\{u_x^n\}$ and
$\{\rho^n\}$ are uniformly bounded on $(0,T)\times(a,b)$. Using
Lemma 4.6, there exist a subsequences $\{u_x^{n_k},
\rho^{n_k}\}_{k=1}^{\infty}\subset\{u_x^n,\rho^n\}_{n=1}^{\infty}$
 for a.e. $(t,x)\in(0,T)\times(a,b)$ and two Borel probability measures $\mu_{(t,x)},\nu_{(t,x)}$
 on $\mathbb{R}$ such that for each $F\in
C(\mathbb{R})$ we have
$$F(u^{n_k}_{x})\rightharpoonup^*\overline{F}_1, \ \ F(\rho^{n_k})\rightharpoonup^*\overline{F}_2 \ \ \text{in}\ \
L^\infty((0,T)\times(a,b)),$$ where
$\overline{F}_1(t,x)=\int_{\mathbb{R}}F(y)d\mu_{(t,x)}$ and
$\overline{F}_2(t,x)=\int_{\mathbb{R}}F(y)d\nu_{(t,x)}$ for a.e.
$(t,x)\in(0,T)\times(a,b)$. By Lemma 4.5, we have that
$\mu_{(t,x)}=\delta_{u_x(t,x)}$ and $\nu_{(t,x)}=\delta_{\rho(t,x)}$
for a.e. $(t,x)\in(0,T)\times(a,b).$ Then, in view of Lemma 4.6
 and from the arbitrariness of $T,a,b$, we obtain
\begin{align}
u_x^{n_k}\rightarrow u_x,\ \
\rho^{n_k}\rightarrow \rho\ \ \ \text{in} \ \
L^2_{loc}(\mathbb{R}_+\times\mathbb{R}).
\end{align}

\par
With the above preparations, we now conclude the proof of the
Theorem 1.1. Let $u$ be the limit of the  approximate solutions
$u^{n_k}$ as $n_{k}\rightarrow\infty$. It then follows from Lemma
2.2, Lemma 4.1 and Remark 4.1 that $u\in
C([0,+\infty)\times\mathbb{R})\cap
L_{loc}^{\infty}(\mathbb{R_+},H^1(\mathbb{S}))$ and $\rho\in
L^{\infty}(\mathbb{R_+},L^2(\mathbb{S}))$ hold. From (4.20), we have
that
$$\partial_{x}g*\left(2\mu_{0}u^{n}+\frac{1}{2}(u^{n}_{x})^{2}+\frac{1}{2}(\rho^{n})^{2}\right)\rightarrow
\partial_{x}g*(2\mu_{0}u+\frac{1}{2}u_{x}^{2}+\frac{1}{2}\rho^{2})$$ in the
sense of distributions on $\mathbb{R_+}\times\mathbb{R}$. This shows
that $z$ satisfies (2.5) in the sense of distributions on
$\mathbb{R_+}\times\mathbb{R}$.

From Lemma 4.2, we have $z\in
C^w_{loc}(\mathbb{R}_+;H^1(\mathbb{S})\times L^2(\mathbb{S}))$.
Consequently, we will prove that $z\in
C(\mathbb{R}_+;H^1(\mathbb{S})\times L^2(\mathbb{S}))$. Note that
$u\in C((0,\infty); L^{\infty}(\mathbb{S})),$ it is enough to show
that $\int_{\mathbb{S}}(u_{x}^{2}+\rho^{2})dx$ is conserved in time.
Indeed, if this holds, then
\begin{align*}
&\|z(t)-z(s)\|^2_{H^1(\mathbb{S})\times L^2(\mathbb{S})}\\= \
&\|u(t)-u(s)\|_{L^{2}(\mathbb{S})}^{2}+\|u_{x}(t)-u_{x}(s)\|_{L^{2}(\mathbb{S})}^{2}
+\|\rho(t)-\rho(s)\|_{L^{2}(\mathbb{S})}^{2}
\\= \ &\|u(t)-u(s)\|_{L^{2}(\mathbb{S})}^{2}+\|u_{x}(t)\|_{L^{2}(\mathbb{S})}^{2}+\|u_{x}(s)\|_{L^{2}(\mathbb{S})}^{2}
-2(u_{x}(s), u_{x}(t))_{L^{2}(\mathbb{S})}\\
+&\|\rho(t)\|_{L^{2}(\mathbb{S})}^{2}+\|\rho(s)\|_{L^{2}(\mathbb{S})}^{2}
-2(\rho(s), \rho(t))_{L^{2}(\mathbb{S})}\\
= \ &
\|u(t)-u(s)\|_{L^{2}(\mathbb{S})}^{2}+2(\|u_{0,x}\|_{L^{2}(\mathbb{S})}^{2}+\|\rho_{0}\|_{L^{2}(\mathbb{S})}^{2})
-2((u_{x}(s), u_{x}(t))_{L^{2}(\mathbb{S})}+(\rho(s),
\rho(t))_{L^{2}(\mathbb{S})}),
\end{align*}
$\forall \ t,s\in\mathbb{R}_+.$ Since
$\|u(t)-u(s)\|_{L^{2}(\mathbb{S})}^{2}\rightarrow 0$ and
$$(u_{x}(s), u_{x}(t))_{L^{2}(\mathbb{S})}+(\rho(s),
\rho(t))_{L^{2}(\mathbb{S})}\rightarrow
\|u_{x}(t)\|_{L^{2}(\mathbb{S})}^{2}+\|\rho(t)\|_{L^{2}(\mathbb{S})}^{2}
=\|u_{0,x}\|_{L^{2}(\mathbb{S})}^{2}+\|\rho_{0}\|_{L^{2}(\mathbb{S})}^{2},$$
as $s\rightarrow t,$ we have $z\in
C(\mathbb{R}_+;H^1(\mathbb{S})\times L^2(\mathbb{S}))$.

The conservation of $\int_{\mathbb{S}}(u_{x}^{2}+\rho^{2})dx$ in
time is proved by a regularization technique. Denote
$f_n=f\ast\phi_n.$ By Lemma 4.5 and (4.15)-(4.16), we have
\begin{equation}
\frac{\partial u_{n,x}}{\partial t}-(u+\gamma)\frac{\partial
u_{n,x}}{\partial x}=(\frac{1}{2}u^2_x-2\mu_{0}u-\frac{1}{2}\rho^{2}
+2\mu_{0}^{2}+\frac{1}{2}\mu_{1}^{2})\ast\phi_{n}+\tau_{n}
\end{equation} and
\begin{equation}
\frac{\partial\rho_{n}}{\partial
t}-(u+\gamma)\frac{\partial\rho_{n}}{\partial
x}=(u_x\rho)\ast\phi_{n}+\sigma_{n},
\end{equation}
Multiplying (4.21) with $ u_{n,x}$, we obtain by integration
\begin{align}
&\frac{1}{2}\frac{d}{dt}\int_{\mathbb{S}}u_{n,x}^{2}dx
-\frac{1}{2}\int_{\mathbb{S}}(u+\gamma)(u_{n,x}^{2})_{x}dx\\
\nonumber= \
&\int_{\mathbb{S}}u_{n,x}(\frac{1}{2}u^2_x-2\mu_{0}u-\frac{1}{2}\rho^{2}
+2\mu_{0}^{2}+\frac{1}{2}\mu_{1}^{2})\ast\phi_{n}dx+\int_{\mathbb{S}}u_{n,x}\cdot\tau_{n}dx.
\end{align}
Multiplying (4.22) with $\rho_{n}$, we get by integration
\begin{equation}
\frac{1}{2}\frac{d}{dt}\int_{\mathbb{S}}\rho_{n}^{2}dx
-\frac{1}{2}\int_{\mathbb{S}}(u+\gamma)(\rho_{n}^{2})_{x}dx=\int_{\mathbb{S}}\rho_{n}(u_x\rho)\ast\phi_{n}dx+
\int_{\mathbb{S}}\rho_{n}\sigma_{n}dx.
\end{equation}
Adding (4.23) and (4.24), we have
\begin{align*}
&\frac{1}{2}\frac{d}{dt}\int_{\mathbb{S}}(u_{n,x}^{2}+\rho_{n}^{2})dx\\
= \
&\frac{1}{2}\int_{\mathbb{S}}(u+\gamma)(u_{n,x}^{2})_{x}dx+\int_{\mathbb{S}}u_{n,x}(\frac{1}{2}u^2_x-2\mu_{0}u-\frac{1}{2}\rho^{2}
+2\mu_{0}^{2}+\frac{1}{2}\mu_{1}^{2})\ast\phi_{n}dx+\int_{\mathbb{S}}u_{n,x}\cdot\tau_{n}dx\\
+ \
&\frac{1}{2}\int_{\mathbb{S}}(u+\gamma)(\rho_{n}^{2})_{x}dx+\int_{\mathbb{S}}\rho_{n}(u_x\rho)\ast\phi_{n}dx+
\int_{\mathbb{S}}\rho_{n}\sigma_{n}dx.
\end{align*}
As for fixed $T>0$, $u$, $u_x$ and $\rho$ are bounded in
$(0,T)\times \mathbb{S}$. Let $n\rightarrow\infty$, on account of
Lebesgue's dominated convergence theorem, we get
$$\frac{d}{dt}\int_{\mathbb{S}}(u_{x}^{2}+\rho^{2})dx=0.$$
This completes the proof of Theorem 1.1.

\bigskip

\noindent\textbf{Acknowledgments} This work was partially supported
by NNSFC (No. 10971235), RFDP (No. 200805580014), NCET-08-0579 and
the key project of Sun Yat-sen University.

\end{document}